\documentclass[11pt,leqno]{amsart}

\makeatletter
\@namedef{subjclassname@2020}{\textup{2020} Mathematics Subject Classification}
\makeatother

\usepackage[top=2.5 cm,bottom=2 cm,left=2.5 cm,right=2 cm]{geometry}
\usepackage[utf8]{inputenc}
  \usepackage{color}

\usepackage{esint}
\usepackage{graphicx}

\usepackage{hyperref}
\usepackage{color}
\usepackage{epsfig}
\usepackage{ulem}
\usepackage{lettrine}
\usepackage{verbatim}
\usepackage{array}
\usepackage{varwidth}
\usepackage{mathptmx} 
\usepackage{setspace}
\usepackage[numbers,compress]{natbib}
\usepackage[dvipsnames]{xcolor}
\usepackage{bm}
\usepackage{comment}
\usepackage{physics}
\numberwithin{equation}{section}

\newtheorem*{maintheorem*}{Main Theorem}
\allowdisplaybreaks
\numberwithin{equation}{section}

\renewcommand{\i}{\ifmmode\mathit{\mathchar"7010 }\else\char"10 \fi}
\renewcommand{\j}{\ifmmode\mathit{\mathchar"7011 }\else\char"11 \fi}
\newcommand{\R}{\mathbb{R}}

\newcommand{\p}{\partial}

\def\begi{\begin{itemize}}
\def\endi{\end{itemize}}
\def\bega{\begin{array}}
\def\enda{\end{array}}

\newcommand{\xxi}{\bm{\xi}}
\newcommand{\eeta}{\bm{\eta}}

\def\u{\mathbf{u}}

\def\b{\mathbf{b}}
\def\x{\mathbf{x}}
\def\f{\mathbf{f}}
\def\y{\mathbf{y}}

{%

\begin{enumerate}}%
{\end{enumerate}}

{%

\begin{enumerate}}%
{\end{enumerate}}

\begin{document}

\title[Bond-based peridynamics in fluids]{Bond-based peridynamics, a survey prospecting nonlocal theories of fluid-dynamics}

\author[N. Dimola]{Nunzio Dimola}
\author[A. Coclite]{Alessandro Coclite}
\author[G. Fanizza]{Giuseppe Fanizza}
\author[T. Politi]{Tiziano Politi}

\address[Nunzio Dimola, Alessandro Coclite, Tiziano Politi]{\newline
Dipartimento di Ingegneria Elettrica e dell'Informazione, \newline Politecnico di Bari, Via Re David 200 -- 70125 Bari, Italy}
\email[]{nunzio.dimola@poliba.it, alessandro.coclite@poliba.it, tiziano.politi@poliba.it}

\address[Giuseppe Fanizza]{\newline Instituto de Astrofisíca e Ci\^encias do Espa\c{c}o, Faculdade de Ci\^encias, \newline Universidade de Lisboa, Edificio C8, Campo Grande, P-1749-016, Lisbon, Portugal}
\email[]{gfanizza@fc.ul.pt}

\date{\today}



\begin{abstract}
Peridynamic (PD) theories have gained widespread diffusion among various research areas, due to the ability of modeling discontinuity formation and evolution in materials. Bond-Based Peridynamics (BB-PD), notwithstanding some modeling limitations, is widely employed in numerical simulations, due to its easy implementation combined with physical intuitiveness and stability. In the present paper, several aspects of bond-based peridynamic models have been reviewed and investigated. A detailed description of peridynamics theory, applications, and numerical models has been presented. Employed BB-PD integral kernels have been displayed together with their differences and commonalities; then, some consequences of their mathematical structure have been discussed. The kinematic role of nonlocality, the relation between kernel structure and material impenetrability, and the role of  PD kernel nonlinearity in crack formation prediction have been critically analyzed and commented on. Finally, the idea of extending BB-PD to fluids is proposed and presented in the framework of fading memory material, drawing some perspectives for a deeper and more comprehensive understanding of the peridynamics in fluids.
\end{abstract}

\maketitle

\section{Introduction}
\label{sec:intro}

\begin{figure}
\centering
\includegraphics[scale=0.225]{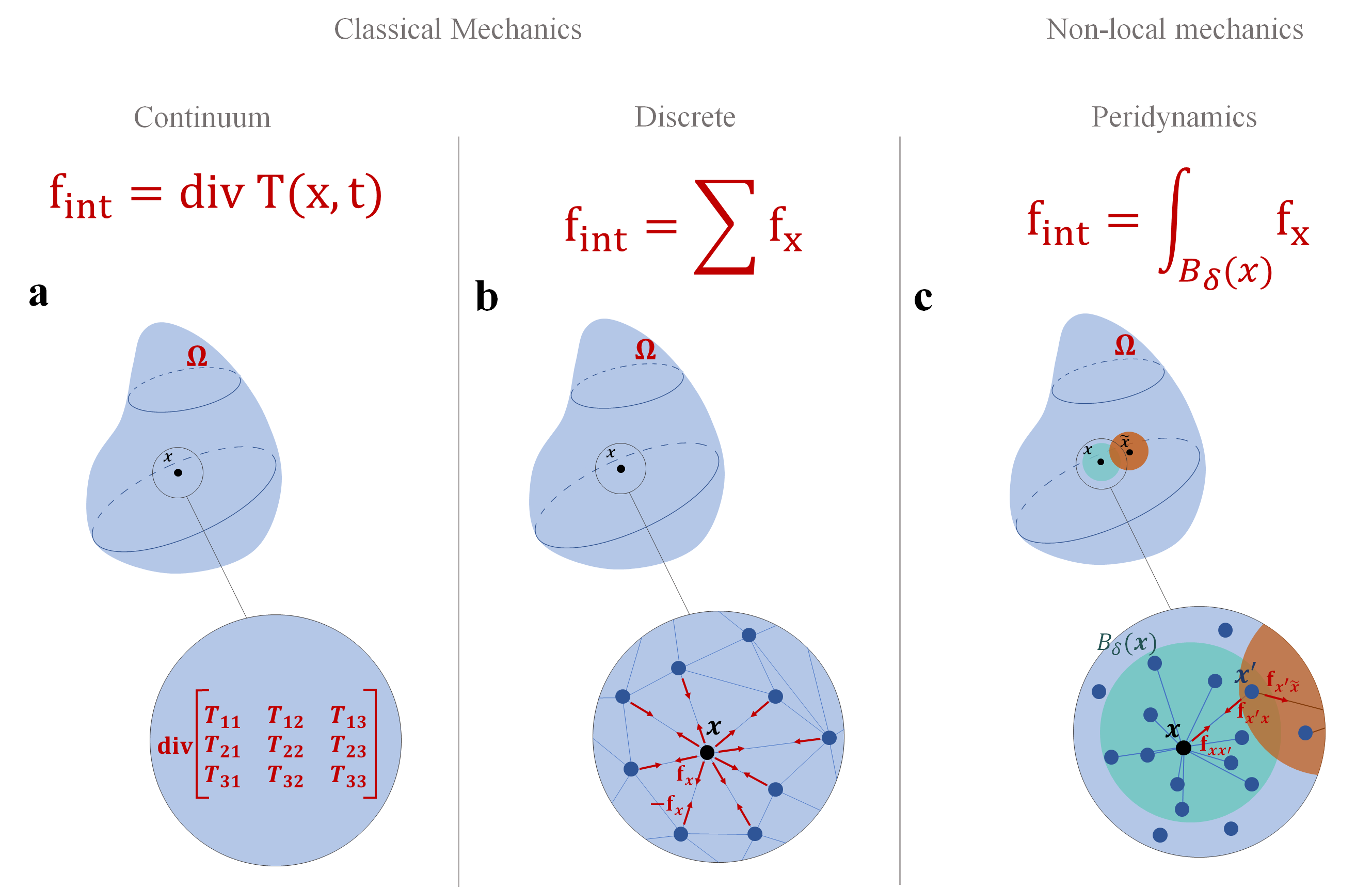}
\caption{{\bf Descriptions of material's points dynamics in continuum and discrete mechanics versus peridynamics.} {\bf a.} Pictorial definition of internal stresses in continuum mechanics hypothesis. ${\bf f}_{int}$ is obtained by $\textrm{div}\bm{T}(\x,t)$, being $ \bm{T}(\x,t) $ the stress tensor. {\bf b.} The material body is discretized by a network of edges connecting adjacent points. The summation of all local forces acting on $\x$ return the local value of $ \f_{int} $. {\bf c.} Within the peridynamics hypothesis, every point $ \x $ is linked with all the points falling in its finite-range horizon, and $ \f_{int} $ is given by the integral of nonlocal forces acting on $\x$.}
\label{confronto}
\end{figure}

Peridynamics (PD) is a nonlocal continuum mechanics theory introduced by Silling \cite{Sill} to mathematically describe fracture formation and development in elastic materials. To overcome limitations imposed by the classical theory, which results inadequate in treatment of spatial discontinuities that eventually occur in material bodies, in the peridynamic equation of motion, an integral operator takes the place of spatial derivatives. The process of integration holds its consistency even with very irregular functions, so that discontinuous displacements (i.e. cracks) are still significant in the peridynamics framework. The growing interest of the scientific community in peridynamics \cite{javili2019peridynamics} is motivated by the fact that such a theory could constitute a linking connection between atomistic theories of matter and classical continuum mechanics. Consequently, micro-scale phenomena, like wave dispersion~\cite{CDFMRV,seleson2009peridynamics,butt2017wave,bavzant2016wave}, cracks (e.g. \cite{ha2010studies,agwai2011predicting,ni2019static,lipton2014dynamic,silling2010crack}), intra-granular fracture \cite{behzadinasab2018peridynamics}, etc. could be modeled, even in macroscopic structures, by a suitable tuning of the \textit{peridynamic horizon} length, which rules the extent of nonlocal interactions between a point of the body and the surrounding ones \cite{askari2008peridynamics}. 
Firstly, the so-called \textit{bond-based} (BB) peridynamic has been introduced in \cite{Sill}, in which internal forces between a point and all the other ones inside its peridynamic horizon are modeled as a central forces field. We can figure such a forces field as a network of bonds linking each point of the body with every point within its horizon. Then, to overcome bond-based PD modeling limitations for the Poisson’s ratio \cite{madenci2014peridynamic,macek2007peridynamics,sarego2016linearized,zaccariotto2015examples} ($1/3$ for plane stress and $1/4$ for plane strain in 2D bodies, 1/4 for 3D ones), a generalized peridynamics theory, \textit{state-base} peridynamics, has been formulated in \cite{silling2007peridynamic}. In this context, the force exchanged between a point and another one in its horizon does not depend solely on their bond extension, but also on the deformation state of all other bonds relative to the horizon.
It should be noted that peridynamics theories differ from classical mechanics for the presence of such finite-range bonds between any two points of the material body, which is a feature that assimilates such formulations to discrete mesoscale theories of matter. In the PD framework, physical bodies are considered as formed by a continuous network of points exchanging momentum within a fixed interaction distance $\delta > 0$: the horizon radius; this change of paradigm, much closer to molecular dynamics than macroscopic bodies one, allows to abandon the local concept of stress tensor and move on to the concept of \textit{pairwise force} (see Fig.\ref{confronto} \textbf{a, c}). In a Lagrangian framework, best suited for finite deformations, the peridynamic horizon is fixed in the reference configuration and deforms with the body \cite{madenci2014peridynamic}. Conversely, when viscoelastic materials or fluids (which are naturally subjected to very large deformations), are considered, it is more physically reasonable to maintain the PD horizon fixed in time, so that its shape never changes as the body deforms. This approach is known in the literature as Eulerian \cite{silling2017modeling} or Semi-Lagrangian \cite{behzadinasab2020semi} peridynamics.

Due to its flexibility, beyond already mentioned research fields (crack formation and propagation in elastic material, wave dispersion in material, intra-granular fracture), the peridynamics nonlocal approach to discontinuities has found applications in several research areas. In geomechanics, PD has been employed in water-induced soil cracks \cite{ni2020hybrid,zhou2020hydromechanical}, geomaterial failure \cite{song2019peridynamics}, rocks fragmentation \cite{panchadhara2017modeling}, \textit{et cetera} (see \cite{zhou2021state}). In biology, long-range interactions in living tissues \cite{lejeune2017modeling}, cellular ruptures, cracking of bio-membranes \cite{taylor2016peridynamic}, and so on \cite{javili2019peridynamics}, have been analyzed within peridynamics framework. In \cite{bobaru2010peridynamic,bobaru2012peridynamic, oterkus2014peridynamic} a peridynamic theory for thermal diffusion has been introduced in order to model heat conduction in materials featuring spatial discontinuities, defects, inhomogeneities, and cracks. Advection-diffusion behavior in PD has been studied in \cite{foster2019nonlocal} for modeling fingering phenomenon for multiphase fluids in porous media, and in \cite{zhao2018construction}, a constructive PD model for transient advection-diffusion problems has been presented.  Being a plethora of physical phenomena modeled consistently by peridynamics, various multi-physics analysis have been performed within PD scheme, e.g  micro-structural analysis \cite{buryachenko2020generalized}, fatigue \cite{hu2017peridynamics} and heat conduction \cite{oterkus2012peridynamic}in composite materials, galvanic corrosion in metals \cite{zhao2021peridynamic}, electricity-induced cracks in dielectric materials \cite{wildman2015dynamic}, and so forth.

Significant similarities between bond-based peridynamics theory and discrete theory of materials can be pointed out. Indeed, in discrete theories of matter, a material body is thought as formed by a discrete network of points linked spring-like bonds. Each point exert long-range elastic interaction with all others with which shares a bond. Thus, each point belonging to this discrete network can be associated with a discrete set of finite-range interacting points, equivalent to the PD horizon (see Fig.\ref{confronto} \textbf{b}). Both in PD and discrete theories, linear and angular momentum balance is not guaranteed locally, so a process of summation, for discrete network, or integration, for continuous one, over all points that constitute the material body is necessary. Given such similarities between material PD theory and discrete one, it is not surprising that most of the scientific efforts have been directed toward a numerical implementation of the peridynamic equation. Apart from PD theory, discontinuities formation and development in material structures have been more or less satisfactory addressed by enhanced versions of classical continuum mechanics numerical methods (like Smoothed Particle Hydrodynamics — SPH  \cite{randles1996smoothed,ren2021improved}, extended Finite Element — xFEM \cite{moes1999finite}, Cohesive Zone Methods — CZM \cite{rocha2020numerical}, and so on); but, peridynamics appeal stems from the lack of external criteria and adjoined degrees of freedom to the numerical schemes. Neglecting some exceptions, PD numerical analysis can be grossly divided into two macro-area: PD finite element analysis and mesh-free methods or, equivalently, quadrature methods \cite{javili2019peridynamics,shojaei2022hybrid}. Such a division is solely for summarizing purposes, and a plethora of subdivisions, parallelisms, and couplings can be found in the literature. Moreover, further numerical methods, like spectral methods \cite{lopez2021space,CFLMP,jafarzadeh2020efficient}, \textcolor{blue}{considering or neglecting the volume penalization step~\cite{LP,lopez2022non}}, Boundary Element Methods (BEM) \cite{liang2021boundary}, have been developed recently in the context of peridynamics, enlarging the range of available PD numerical tools.   In~\cite{silling2005meshfree}, a mesh-free method has been developed to numerically solve the PD equation for so-called prototype micro elastic brittle (PMB) materials, showing that the stability criterion for the proposed numerical scheme weakly depends on space discretization, but results principally dictated by peridynamic horizon size. In \cite{emmrich2007analysis}, mathematical well-posedness of the elastic one-dimensional PD Cauchy problem has been proved and a proposed numerical quadrature method for the PD integrodifferential equation (IDE) was shown to converge to analytical solutions, both for continuous and discontinuous initial conditions. Convergence in mesh-free peridynamics simulations has been analyzed in \citep{seleson2016convergence} for 1D, 2D bond-based PD, and 3D state-based PD, considering static problems; convergence of the numerical scheme, under mesh refinement ($\delta$ kept constant), to a manufactured nonlocal analytical solution has been shown. Two methods for mitigating the so-called \textit{partial volume} effect, i.e. the intersection discrepancy between the circular/spherical shape of the PD horizon and the polygonal shape of grid's cells have been proposed: one based on more accurate analytic calculations of such intersection's contribution, the other based on employing compactly supported smooth micro-modulus functions, i.e functions vanishing at the horizon boundary, so that, the intersected volumes contribution is minimized. In \cite{bessa2014meshfree}, close connections between peridynamics and classical mesh-free methods have been established. State-based PD is shown to be a numerically faster case of the so-called Reproducing Kernel Particle Methods (RKPM), so that, the same approach employed in RKPM boundaries treating could be extended to state-based PD, developing the \textit{reproducing kernel} peridynamics. Mesh refinement strategy for numerical analysis of 1-D peridynamic elastic bodies has been introduced by \cite{bobaru2009convergence}; convergence of the numerical solution to the classical elastic one has been proven for vanishing PD horizon ($\delta \rightarrow 0$) and compared for various employed micro-modulus function shapes; so-called \textit{visibility criterion} was introduced as refinement approach. 

It is a well-known fact that, in correspondence to the interface between refined and coarsen mesh regions, spurious phenomena like volume losses, ghost forces and surface softening may occur \cite{le2018surface,bobaru2011adaptive, dipasquale2014crack}. In \cite{ren2016dual} the concept of \textit{dual horizon} has been firstly introduced: the total internal force acting at a point of the body is split between a reaction force component, due to the point's interactions with its neighborhood (because of Newton's Third Law), and an active force component exerted on the point by its dual horizon. In such a way, peridynamics forces are defined consistently also in correspondence with refined mesh interfaces, where different horizon length coexists. The dual horizon concept has been applied in \cite{gu2017voronoi} to discretized bond-based PD, together with Voronoi-based adaptive mesh refinement, showing good agreement with standard peridynamics test cases. In \cite{shojaei2018adaptive} a careful approach is adopted in mesh adaptive refinement by the introduction of fictitious nodes near the interface region between the coarse and refined mesh. As a result, spurious phenomena are considerably reduced. In \cite{henke2014mesh}, mesh influence on PD numerical solution has been explored in bond-based PD framework. The onset of periodic patterns in the integration error for 1D configuration, accuracy loss caused by perturbations of quadrature points location in 2D configuration, and unrealistic crack patterns in 3D impact simulations have been observed and directly linked to grid Cartesian structure and simplistic quadrature strategies, pointing out that accurate discretization schemes mush be employed for preventing mesh sensitivity in PD numerical simulations.

Concurrently, Finite Element Methods (FEM) have been extended to peridynamics theory. In \cite{macek2007peridynamics}, peridynamics extensibility to finite elements method coding has been proved and tested, using ABAQUS\textsuperscript{\textregistered} FEM code, for impact simulations, showing great prediction capabilities. Given the wide diffusion of commercial FEM codes in engineering, the interest in couplings between FEM code and PD theory has rapidly increased. In \cite{kilic2010coupling}, a coupling between FEM and peridynamics has been carried out, to take advantage of the peculiarities of both FEM and PD approach, i.e. fast numerical efficiency and inherent crack prediction capabilities, respectively. The Peridynamics approach has been employed only for domain portions where damage was expected to occur, while the FEM approach has been employed for no-expected damage zones. Application of continuous and discontinuous Galerkin finite element methods to PD has been explored in \cite{chen2011continuous} and validated against 1D peridynamics exact solutions. In \cite{liu2017hp}, a fast and cost-efficient Galerkin method developed in \cite{wang2012fast} has been extended and improved for 1D linearized PD static problems providing \textit{hp}-adaptivity, reducing computational efforts and memory usage. In \cite{huang2019finite}, the implementation of a coupled PD-FEM approach in commercial software ABAQUS\textsuperscript{\textregistered} has been performed, employing mesh coupling techniques developed in  \cite{zaccariotto2017enhanced,zaccariotto2018coupling,galvanetto2016effective}, with good results; while, recently, in \cite{zhang2022coupled} a PD-FEM coupled approach has been carried out with the commercial software ANSYS\textsuperscript{\textregistered} for fatigue prediction in materials.

Although bond-based PD is the most dated formulation of peridynamics and suffers of the discussed Poisson's ratio limitation, it is yet employed in several numerical studies on peridynamics. Physical intuitiveness, less computational efforts and crack simulation stability concerning state-based PD, makes the choice of bond-based integral kernel still attractive \cite{zheng2020bond,han2019review}, so much that several extension of bond-based PD has been proposed (for a review refer to \cite{han2019review}).

In bond-based numerical analysis, different types of integral kernels have been employed to model different kinds of nonlocal interactions, characteristics of the analyzed material, and/or structure. Various micro-modulus functions have been studied and tested in the classical or \textit{linearized} version of BB-PD, resulting in a great variety of constitutive models. At the same time, the fundamental structure of the BB peridynamics kernel, with few exceptions \cite{nanoBobaru,lopez2021space,CFLMP}, has not varied so much  from the prototype micro elastic brittle (PMB) material (e.g. \cite{chen2019influence,ha2010studies,kilic2008peridynamic} among others), which is fundamentally a \textit{linear} model. A clarification must be done: with a linear (nonlinear) model we design all those PD kernels which are linear (nonlinear) in the displacement variable and undergo a \textit{finite} displacement; with \textit{linearized} model we mean every BB-PD integral kernel linearized (via Taylor expansion) with respect to the relative displacement variable, i.e. that undergo, for small relative displacements, a first-order approximation of the pairwise force.

The great flexibility in choosing a peridynamics kernel, however, addresses a fundamental question: how do physical principles reflect on the structure of the peridynamic equations? Or, conversely, which physic emerges according to the different kernels structures with which peridynamic theory is presented?
This paper aims to shed light on the questions above to highlight common features between widespread linear and nonlinear peridynamics bond-based models for hyperelastic materials and to shape some future perspectives for a wider understanding of kernel's structure-related physical and numerical peculiarities.

Basic ideas, concepts, and physical assumptions of peridynamics bond-based theory in Lagrangian formalism are illustrated in Section \ref{sec:foundation}. A review of bond-based integral kernels employed in the literature is presented in Section \ref{sec:PD models}, underlying their features, structures, and properties. In Section \ref{subsec:non-loc}, the kinematics aspects of peridynamics nonlocality property are highlighted, compared with those of classical continuum mechanics, and related to the concept of the finiteness of bond stretch. The influence of PD kernel structure on the basic requirement of impenetrability of the matter is explored in Section \ref{subsec:non-pen}. The potential role of peridynamic kernel nonlinearity as a substitute for the classical critical-stretch criterion in material damage prediction is formulated in Section \ref{subsec: nonlincrak}. Finally, in Section \ref{sec: perifluid}, some perspectives are traced for the application of bond-based Peridynamic to fluids. Fluids and solids are compared in their constitutive behavior, then, by making use of the merging concept of fading memory, some aspects of an extension of PD formalism toward fluids have been formally highlighted.

\section{Peridynamics Foundation}
\label{sec:foundation}
\begin{figure}
\centering
\includegraphics[scale=0.4]{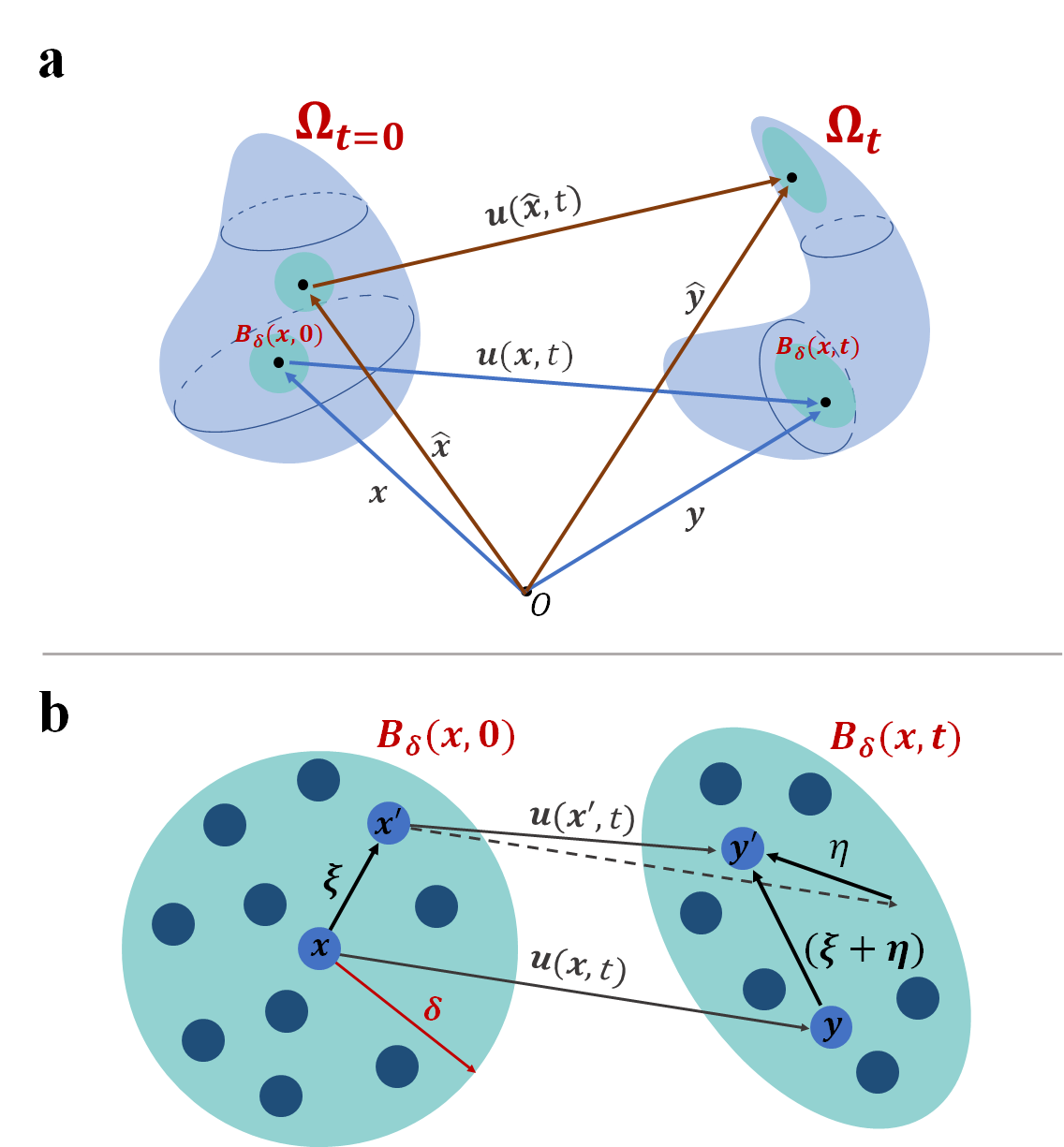}
\caption{{\bf Points kinematic and range of interactions in peridynamics.} {\bf a.} All points in the reference configuration $\Omega_{t=0}\equiv \Omega_{0}$ are labeled through a ray-vector starting in the origin of $\Omega_{0}$. The displacement functions $\u(\x,t)$ associate a point into the deformed configuration $\Omega_t$ to each point of $\Omega_{0}$. Indeed, at $t>0$ a point in $ \Omega_t$ is completely characterized by ${\bf x}(t)=\x(0)+\u(\x,t)$
{\bf b.} Deformation in time of the ball $B_{\delta}(x)$ made up by all points $ \x' $ such that $ \abs{\x-\x'}\leq \delta $.}
\label{cinematic}
\end{figure}

Let us consider a subset $\Omega \subset \R^{n}$ with $n\in \{ 1,2,3 \}$ representing a material body with constant density $\rho>0$. At time $t=0$ such material body encloses the volume $\Omega_{0} \subset \R^{n}$ taken as its reference configuration. Note that, the reference configuration may correspond either to the stress-free configuration or to a given configuration of the body taken as reference. As discussed, in the perydynamics hypothesis, each point of the material body interacts with all of the points $\x'$ falling within a certain neighborhood as such $d(\x,\x')\leq\delta$, being $\delta>0$ and $d(\cdot,\cdot)$ a suitable distance on $\Omega_0$. We refer to this set as $B_\delta(\x)$ (see Fig.\ref{cinematic}b); in the literature it is usually referenced either as the \textit{horizon} (e.g in \cite{Sill,silling2005meshfree,ren2016dual}), or the \textit{family} of $\x$ (e.g. in \cite{madenci2014peridynamic,chen2016constructive,madenci2016peridynamic}).
Here, the kinematics of $\x$ is given in terms of its displacement with respect to the reference position, namely $\u(\x, t): \Omega_{0} \times \mathbb{R}^{+} \rightarrow \mathbb{R}^{n}$. As a consequence, the position of $\x$ at a certain $t$ is given by $\y(\x,t) := \x+\u(\x, t)$ (see Fig.\ref{cinematic}a). Moreover, for each couple of interacting points, the length of the pairwise bond relatively to the initial configuration is tracked in time by the relative strain $s(\x,\x',t)$, reading:
\begin{equation}
s\left(\x, \x', t\right)=\frac{\left|\u\left(\x^{\prime}, t\right)-\u(\x, t)\right|}{\left|\x^{\prime}-\x\right|},
\end{equation}
where $ |\cdot| $ is the Euclidean norm. Such pairwise bonds length varies in time responding to the force per unit volume squared $\f(\x',\x,\u'(\x),\u(\x),t)$. Note that the dependence of $\u$ on $t$ has been omitted to simplify the notation. Alongside $\f$, an external forcing term is provided, $\mathbf{b}(\x,t)$, so that the Lagrangian balance equation for the momentum of $ \x $ reads:
\begin{equation}
{\rho \u_{tt}(\x,t) = {\bf F}(\x,t)}\, .
\end{equation}
${\bf F}(\x,t)$ is the sum of all of the internal and external per-unit-volume forces acting on $\x$:    
\begin{equation}
{{\bf F}(\x, t)=\int_{\Omega_0 \cap B_{\delta}(\x)} \f\left( \x', \x, \u \left( \x' \right),\u(\x) \right) dV_{\x'}+\b(\x, t)}\, .
\end{equation}
Within the hypotheses of \textit{homogeneity} of the material and \textit{invariance with respect to rigid motion}, the pairwise interaction reduces to a function of $\x'-\x$ and $\u(\x')-\u(x)$. This momentum balance equation is enriched by initial conditions and the Cauchy problem responding to nonlocal interactions reads:
\begin{equation}
\begin{cases}
\rho \u_{tt}(\x,t)=({\bf K}\,\u(\cdot,t))(\x)+\b(\x,t)\, , &\quad t >0\, , \x \in \Omega_0\, ,\\
\u(\x,0)=\u_0(\x)\, , &\\
\partial_t \u(\x,0)=\mathbf{v}_0(\x)\, , &
\end{cases}
\label{evolution}
\end{equation}
where
\begin{equation}
{({\bf K}\,\u)(\x):= \int_{\Omega_{0}\cap B_\delta(\x)}\mathbf{f}(\x'-\x,\u(\x')-\u(\x))\,dV_{\x'},\quad \x\in\Omega_0}\, .
\end{equation}

Only closed manifolds will be considered, so that no boundary conditions are required. Looking at (\ref{evolution}), it is noteworthy that no spatial derivatives are involved in the evolution equation of $ \x $, so the integrodifferential equation still remains valid for regions characterized by discontinuous displacements (crack, phase changing region, etc.): this represents a great advantage of the peridynamic theory of continuum. 
  
If external forces in Eq.\eqref{evolution} are neglected, it becomes manifest that only the kernel structure of $({\bf K}\,\u)$ influences the whole problem from both mathematical and physical points of view. That is, the mathematical structure of the evolution problem, its physical coherence, and the constitutive properties of the material under consideration derive directly from the integral kernel. In the following, the properties of the integral kernel $ \f(\x'-\x,\u(\x')-\u(x)) $ will be investigated and its relation with the behavior of the material will be deciphered.
Peridynamics integral kernel is shaped up by the nature of internal forces and, indeed, some physical constraints are required. With the intent to ease notation, the following quantities are defined (see Fig.\ref{cinematic}b): ${\bf \xxi} := \x'-\x$ and $\eeta:=\u(\x')-\u(\x)$ so that ${\bf f}(\x'-\x,\u(\x')-\u(\x)) \equiv \f(\xxi,\eeta)$. 

\begin{itemize}
\item \textit{Actio et reactio principle}:
\begin{equation}
{\forall \x \, , \x' \in B_{\delta}(\x): \f(-\xxi,-\eeta)=-\f(\xxi,\eeta)}\, .
\label{linbal}
\end{equation}
The actio et reactio principle, i.e. the \textit{Third Law of Newton}, ensures the conservation of linear momentum of the system composed of mutually interacting particles. 

\item \textit{Angular momentum conservation}: 
\begin{equation}
{\forall \x\, , \x' \in B_{\delta}(\x): (\xxi+\eeta) \times \f(\xxi, \eeta)= {\bf 0}}\, ,
\end{equation}
considering that the relative deformed ray-vector connecting $\x $ and $ \x' $ is $(\xxi+\eeta)$. The condition above is satisfied if and only if the pairwise force density vector has the same direction as the relative deformed ray vector,  
\begin{equation}
\f(\xxi, \eeta)=f(\xxi, \eeta)(\xxi+\eeta)\, , \ \ \ \ \ \   \forall \xxi\, , \eeta \, ,
\label{mombal}
\end{equation}
with $ f(\xxi, \eeta) $ a scalar-valued function.

\item \textit{Hyperelastic material}: hyperelastic is the attribute given to a material such that:
\begin{equation}
\int_{\Gamma} \f(\xxi, \eeta) \cdot d \eeta=0\, , \quad \forall \text{ closed curve } \Gamma, \ \ \ \ \forall\xxi\neq \bm{0}, 
\end{equation}
or, equivalently, by Stokes' theorem
\begin{equation}
\nabla_{\eeta} \times \f(\xxi,\eeta)=\bm{0}\, ,  \ \ \ \ \ \   \forall \xxi\, , \eeta \, ,
\label{rotor}
\end{equation} 
and consequently,
\begin{equation}
\f(\xxi,\eeta)=\nabla_{\eeta}\Phi(\xxi,\eeta)\, , \ \ \ \ \ \   \forall \xxi\, , \eeta \, .
\end{equation} 
In the above equation $ \Phi(\xxi,\eeta) $ is a scalar-valued potential function in $ C^2(\R^n \setminus\bm{\{0\}} \times \R^n) $ \cite{Sill}. Since pairwise force must satisfy the angular momentum conservation, the following condition on the scalar-valued function $ f(\xxi,\eeta) $ is obtained 
\begin{equation}
\frac{\p f(\xxi,\eeta)}{\p \eeta}=g(\xxi,\eeta)(\xxi+\eeta).
\end{equation}
Integrating both sides of the equation, the following condition on $ g(\xxi,\eeta) $ returns~\cite{Sill}
\begin{equation}
\f(\xxi,\eeta)=h(| \xxi+\eeta|,\xxi)(\xxi+\eeta),
\label{hyperel}    
\end{equation}
for $h(| \xxi+\eeta|,\xxi)$ a scalar valued function.
From (\ref{hyperel}), the elastic nature of $ \f $ can be clearly observed, in fact, the interaction force depends only on the initial relative position between points $ \x $ and $ \x' $ and the modulus of their relative position in the deformed configuration $ \Omega_t$ at time $t$,  $ | \xxi+\eeta| $. Under the \textit{isotropy} hypothesis, the general dependence on vector $ \xxi $ can be substituted with a dependence on $ |\xxi| $,
\begin{equation}
    \f(\xxi,\eeta)=h(| \xxi+\eeta|,|\xxi|)(\xxi+\eeta).
\label{hypermod}    
\end{equation}
Bond forces can, thus, be considered as modeling a spring network that connects each point  $ \x \in \Omega_0 $ pairwise with $ \x' \in B_{\delta}(\x) \cap \Omega_0 $. 
\item \textit{Linear elastic material}: if $ |\eeta| \ll 0 $, the peridynamic kernel can be linearized around $ \eeta=\bm{0} $ so that
\begin{equation}
\f(\xxi,\eeta)\approx \f(\xxi,\bm{0})+\left. \frac{\p \f(\xxi,\eeta)}{\p\eeta}\right|_{\eeta=\bm{0}}\eeta;
\end{equation}
a second-order \textit{micro-modulus tensor} can be defined as
\begin{equation}
\bm{C}(\xxi)=\left. \frac{\p \f(\xxi,\eeta)}{\p\eeta}\right|_{\eeta=\bm{0}}=\xxi \otimes \left.\frac{\p f(\xxi,\eeta)}{\p \eeta}\right|_{\eeta=\bm{0}}+f_0,
\end{equation}
where, to simplify notation, $ f_0=f(\xxi,\bm{0}) $. After an application of linear momentum balance (\ref{linbal}), elasticity (\ref{rotor}), and isotropy condition, the micro-modulus tensor can be expressed in the following form \cite{Sill}
\begin{equation}
    \bm{C}(\xxi)=\lambda(|\xxi|)\xxi \otimes \xxi+f_0.
\end{equation}
Thus, for a linearized hyperelastic material, its peridynamic kernel shows the following structure
\begin{equation}
    \f(\xxi,\eeta)= \f(\xxi,\bm{0})+\left(\lambda(|\xxi|)\xxi \otimes \xxi+f_0\right)\eeta.
    \end{equation}

\end{itemize}

\section{Peridynamic models for hyperelastic materials}
\label{sec:PD models}
  Looking at the literature regarding peridynamic elastic theory, different integral kernel (i.e pairwise force $ \f(\xxi,\eeta) $) structures can be found; in the following, a survey of various employed expressions for peridynamic kernels is detailed. For all the listed models, the pairwise force $ \f(\xxi,\eeta) $ vanishes for $ |\xxi| > \delta $; such a condition on $ \f $, even if not explicitly stated, must be considered to hold.   In \cite{Sill}, a simple model for three-dimensional micro-elastic, isotropic, peridynamic material subject to anti-plane shear is proposed to predict the onset of fractures: 
\begin{equation}
\mathbf{f}(\boldsymbol{\eta}, \xxi)=\left\{\begin{array}{ll}
c\;\frac{|\xxi+\eeta|-|\xxi|}{ |\xxi+\eeta|}(\xxi+\boldsymbol{\eta})\, , & \text { if } \; |\xxi+\eeta|-|\xxi| \leq u_{*} \; \text { and } \; |\xxi| \leq \delta \, , \\
\bm{0}\, , & \text{otherwise,} 
\end{array}\right.
\label{mod1}
\end{equation}   
where $ c $ is a constant depending on the material and $ u_* $ is a parameter such that, if the displacement of the relative point in the deformed configuration exceeds the threshold level $ u_* $, the associated bond is broken and fracture can emerge. Another model proposed in \cite{Sill} for isotropic micro-elastic peridynamic material is described by the following potential function
\begin{equation}
    \Phi\left(|\xxi+\eeta|,|\xxi|\right)=\alpha(|\xxi|)\left(|\xxi+\eeta|^2-|\xxi|^2\right)^2,
    \label{pot2}
\end{equation}
so that
\begin{equation}
    \mathbf{f}(\boldsymbol{\eta}, \xxi)=a(|\xxi|)\left(|\xxi+\eeta|^2-|\xxi|^2\right)(\xxi+\eeta),
    \label{mod2}
\end{equation}
where $ a(|\xxi|) $ is a scalar function.

\begin{figure}
\centering
\includegraphics[scale=0.4]{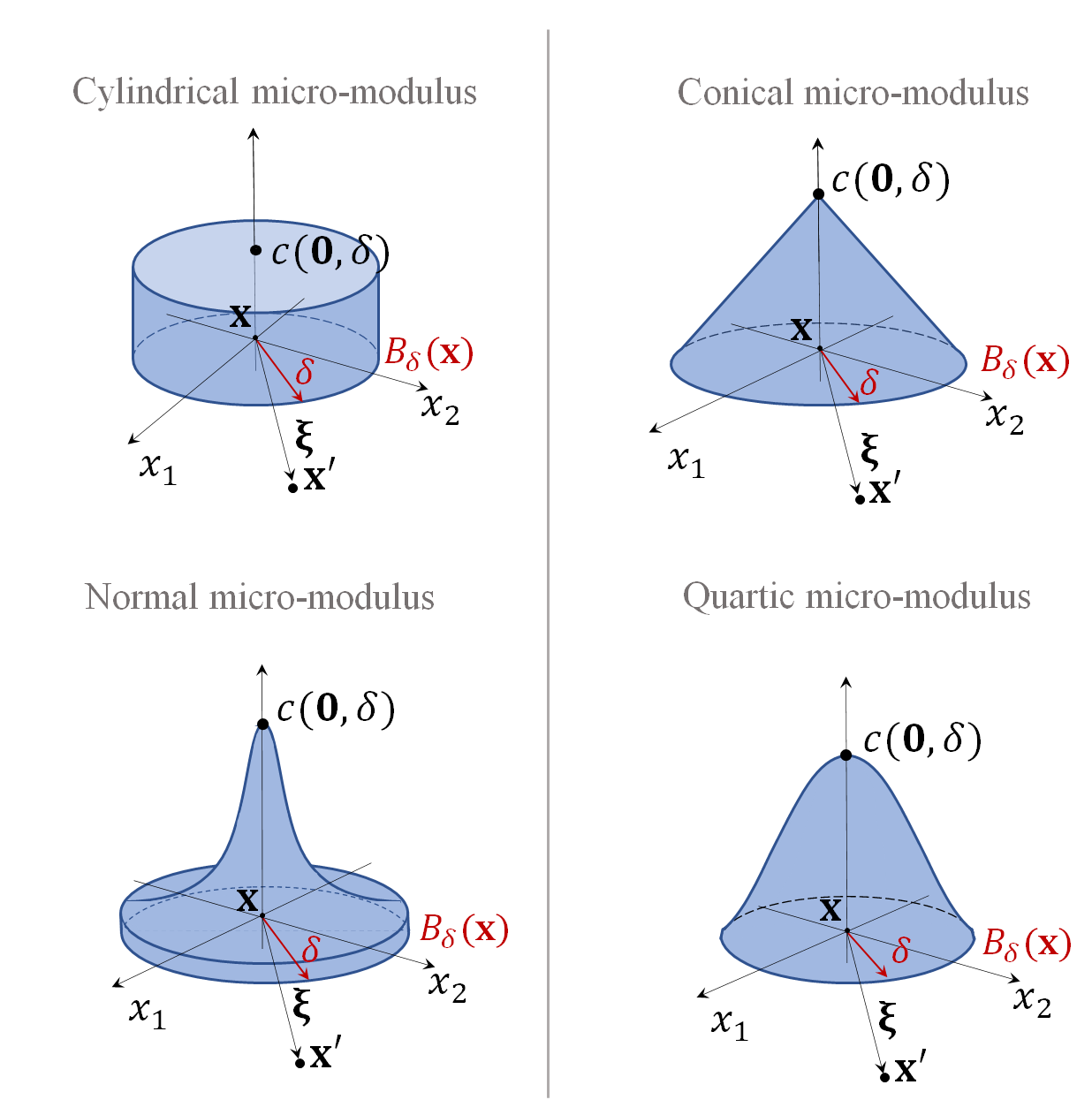}
\caption{{\bf Most common micro-moduli function $c(\xxi,\delta)$; two-dimensional space representation.} \textbf{a)} \textit{Cylindrical} micro-modulus. \textbf{b)} \textit{Triangular} micro-modulus. \textbf{c)} \textit{Normal} micro-modulus \textbf{d)} \textit{Quartic} micro-modulus.}
\label{micromodul}
\end{figure}

In \cite{silling2005meshfree}, a peridynamic kernel for the so-called \textit{prototype micro-elastic brittle} (PMB) material is proposed; the pairwise force for an isotropic PMB material is conceived as linearly proportional to the finite stretch $ s:= (|\xxi+\eeta|-|\xxi|)/|\xxi| $, so that 
\begin{equation}
\mathbf{f}(\boldsymbol{\eta}, \xxi)=f(|\xxi+\eeta|,|\xxi|) \bm{n},
\label{mod3}
\end{equation}   
where $ \bm{n}:=(\xxi+\eeta)/|\xxi + \eeta| $ and where the scalar function $ f $ is defined as follow
\begin{equation}
f=cs\mu(s,t)=c \; \frac{|\xxi+\eeta|-|\xxi|}{|\xxi|}\mu(s,t),
\end{equation}   
with 
\begin{equation}
\mu(s,t)=\left\{\begin{array}{ll}
1\, , & \text { if } s\left(t^{\prime}, \xxi\right)<s_{0}\, , \\
0\, , & \text { otherwise, }
\end{array}\ \ \ \  \text { for all } 0 \leq t^{\prime} \leq t\right.;
\end{equation}
$ c $ is the \textit{micro-modulus constant} and $ \mu(s, t) $ is a function that records if, at a certain time $ t'\leq t $, the bond stretch $ s $ associated to the pair $ (\x,\x') $ has exceeded the critical value $ s_0 $: if so, the bond is considered \textit{broken}, and a zero-valued pairwise force is assigned for all $ t \geq t' $. By comparing the strain energy density value obtained under  isotropic extension respectively with peridynamics and classical continuum theory, the value of $ c $ is readily found \cite{silling2005meshfree}
\begin{equation}
c=\frac{18 k}{\pi \delta^{4}},
\end{equation}
where $k$ is the material bulk modulus. Over time, following the same approach as in \cite{bobaru2010peridynamic} for micro-conductivity function, the micro-modulus constant $ c $ has been generalized to $ c(\xxi,\delta) $ i.e. to a \textit{micro-modulus function}, in order to describe in a more detailed way how pairwise forces intensity distributes over the peridynamic horizon $ B_{\delta}(\x) $; intuitively, forces intensity decrease as the distance between $\x$ and $\x' \in B_{\delta}(\x) $ increases, but it can do it in various ways. The  micro-modulus function is defined as
\begin{equation}
    c(\xxi,\delta):=c(\bm{0},\delta)k(\xxi,\delta)\, ,
\end{equation}
where $ c(\bm{0},\delta) $ is a constant obtained by comparing peridynamic strain density with the classical mechanical theories \cite{chen2019influence}; $ k(\xxi,\delta) $ is a function defined on $ \Omega_0 $ such that (in order to meet momentum conservation and isotropy conditions) \cite{huang2015extended}
\begin{equation}
\left\{\begin{array}{l}
k(\xxi, \delta)=k(-\xxi, \delta)\, , \\
\lim _{\xxi \rightarrow \bm{0}} k(\xxi, \delta)=\max_{\xxi\ \in \R^n}\{ k(\xxi,\delta)\}\, , \\
\lim _{\xxi \rightarrow \delta} k(\xxi, \delta)=0 \, ,\\
\int_{\R^n} \lim _{\delta \rightarrow 0} k(\xxi, \delta) d \x=\int_{\R^n} \Delta(\xxi) d \x=1\, ,
\end{array}\right.
\end{equation}

where $ \Delta(\xxi) $ is the Dirac Delta function. The simplest adopted form for the micro-modulus function is $c(\bm{0},\delta)k(\xxi,\delta)=c\bm{1}_{B_{\delta}(\x')} $, where $ \bm{1}_{A} $: $ X \rightarrow \R $ is the indicator function for the subset $ A \subset X $, defined as
\begin{equation}
\mathbf{1}_{A}(x):= \begin{cases} 1, & x \in A\, , \\ 0, & x \notin A\, , \end{cases}\; \;;
\end{equation}
this kind of micro-modulus function, referred as \textit{cylindrical} (see Fig.\ref{micromodul}a), corresponds to peridynamic kernel in (\ref{mod3}). In \cite{ha2010studies} a  \textit{triangular} micro-modulus function (see Fig.\ref{micromodul}b) is introduced, characterized by $ k(\xxi,\delta) $ to a be a linear function of the following type
\begin{equation}
    k(\xxi,\delta)= \left( 1-\frac{|\xxi|}{\delta} \right)\bm{1}_{B_{\delta}(\x')}. 
\end{equation}
Considered that most common discrete physical systems are characterized by a Maxwell-Boltzmann distribution, in order to include this behavior in peridynamics, a \textit{normal} micro-modulus function (see Fig.\ref{micromodul}c) is proposed in \cite{kilic2008peridynamic} 
\begin{equation}
k(\xxi,\delta)=e^{-(|\xxi| / \delta)^{2}}\bm{1}_{B_{\delta}(\x')};
\end{equation}
while in \cite{huang2015extended}, a \textit{quartic polynomial} micro-modulus (see Fig.\ref{micromodul}d) is proposed 
\begin{equation}
k(\xxi, \delta)=\left(1-\left(\frac{\xi}{\delta}\right)^{2}\right)^{2}\bm{1}_{B_{\delta}(\x')}.
\end{equation}
Therefore, peridynamic kernels like those proposed in \cite{silling2005meshfree}, \cite{huang2015extended}, \cite{ha2010studies} and \cite{kilic2008peridynamic}  can be included in a generalized prototype micro-elastic brittle (PMB) material kernel characterized by the following form
\begin{equation}
    \f(\xxi,\eeta)=c(\xxi,\delta)\underbrace{\frac{|\xxi+\eeta|-|\xxi|}{|\xxi|}}_{s}\mu(s,t);
    \label{brittle}
\end{equation}
an interesting survey and analysis of the effects of different micro-modulus functions on fracture development can be found in \cite{chen2019influence}.
In \cite{nanoBobaru}, peridynamic kernels for failure analysis of material nanostructures are presented: a pairwise force model is proposed for an isotropic Blatz–Ko rubbery membrane and fiber materials characterized by van der Waals interactions, respectively
\begin{equation}
\f(\xxi, \eeta)=\frac{2 c}{|\xxi|}\left(\frac{|\xxi+\eeta|}{|\xxi|}-\left(\frac{|\xxi+\eeta|}{|\xxi|}\right)^{-3}\right) g(|\xxi|) \mu(|\xxi+\eeta|, t)\;\mathbf{n}
\label{nanoplane}
\end{equation}
and
\begin{equation}
\begin{split}
\mathbf{f}(\xxi,\eeta) = \left[\frac{2 c}{|\xxi|}\left(\frac{|\xxi+\eeta|}{|\xxi|} -\left(\frac{|\xxi+\eeta|}{|\xxi|}\right)^{-3}\right)\right] & g(|\xxi|) \mu(|\xxi+\eeta|, t)\\ & -\frac{12 \alpha}{\delta}\left(\frac{\delta}{|\xxi+\eeta|}\right)^{13}
+\frac{6 \beta}{\delta}\left(\frac{\delta}{|\xxi+\eeta|}\right)^{7} \mathbf{n}\, ,
\end{split}
\label{nanofibre}
\end{equation}
where $\alpha, \beta, c $ are constants and $ g(|\xxi|) $is a continuous function that takes into account the possibility of different elastic responses for different characteristic bonds length.
In \cite{chen2016constructive}, the peridynamic elastic material is thought as composed of a superposition of elastic rods with elastic micro-modulus function $ c(\xxi,\delta) $,  that connect pairwise each point belonging to the body. 
Starting from Newton's second laws applied to each rod connecting the point $ \x $ with its family, the following constructive integral kernel is obtained 
\begin{equation}
\mathbf{f}(\boldsymbol{\eta}, \xxi)= c(\xxi,\delta) \frac{|\xxi+\eeta|-|\xxi|}{ |\xxi|^2}\bm{n}.
\label{mod4}
\end{equation}   
In \cite{LP},\cite{lopez2021space}, where one-dimensional and two-dimensional nonlinear elastic materials have been considered in the development of a numerical approximation method for the resolution of the peridynamic problem, the following integral kernel structure of convolution type is proposed \begin{equation}
\f(\xxi,\eeta)= C(\xxi)w(\xxi+\eeta); 
\end{equation}   
the micro-modulus function $ C(\xxi) $ is even while $ w(\xxi+\eeta) $ is odd, both are global Lipschitz continuous functions; specifically, the authors choose $ w(\xxi+\eeta)= (\xxi+\eeta)^{r} $ where $ r >1 $ is an odd number. The peridynamic kernel results as follow
\begin{equation}
\f(\xxi, \eeta)= C(\xxi)(\xxi+\eeta)^r. 
\label{mod5}
\end{equation}  
In \cite{Coclite_2018}, \cite{coclite2021dispersive}, the following nonlinear peridynamic integral kernel in $ \R^n $ for an hyperelastic material is considered by establishing the following nonlinear elastic potential function
\begin{equation}
    \Phi(\xxi,\eeta)=\kappa \frac{|\xxi+\eeta|^p}{|\xxi|^{n+\alpha p}}+\Psi(\xxi,\eeta),
    \label{pot6}
\end{equation}
where $ n $ is the space dimension (practically $ n \in \{1,2,3 \} $), $ p $ and $ \alpha $ real numbers such that
\begin{equation}
    p \geq 2, \quad \alpha \in (0,1),
\end{equation}
and $ \Psi(\xxi,\eeta) $ is a sufficiently smooth perturbation function. Considering that for hyperelastic materials, $ \f=\nabla_{\eeta}\Phi  $ holds, the following structure for the constitutive pairwise force is obtained
\begin{equation}
\mathbf{f}(\boldsymbol{\xi}, \boldsymbol{\eta})=\kappa p \frac{|\boldsymbol{\eta}+\xxi|^{p-2} }{|\boldsymbol{\xi}|^{N+\alpha p}}(\boldsymbol{\eta+\xi})+\psi(\boldsymbol{\xi}, \boldsymbol{\eta}).
\label{mod6}
\end{equation}
It is noteworthy that, for the expression given above, and in general for all pairwise forces defined by an elastic potential function such that
\begin{equation}
{\int_{\mathbb{R}^{N}} \int_{B_{\delta}(\mathbf{0})} \Phi\left(\xxi, \mathbf{u}_{0}(\mathbf{x})-\mathbf{u}_{0}(\x')\right) d \mathbf{x} d \xxi<\infty}\, ,
\end{equation} 
and guaranteed that $\mathbf{u}_{0}\, , \mathbf{v}_{0} \in L^{2}\left(\mathbb{R}^{N} ; \mathbb{R}^{N}\right)$, the associated Cauchy problem in (\ref{evolution}) admits a unique, stable and energy-bounded weak solution, so that the peridynamic equation is well-posed; for a detailed treatment and rigorous proofs of the aforementioned well-posedness of such type of peridynamic models see \cite{Coclite_2018}. 

Then after neglecting minor differences, the pairwise force models listed below could be thought representative of the peridynamic kernel structures commonly adopted in literature for the modeling of hyperelastic materials. An interesting feature of the listed above models is that, even though each of them ensures compliance with the basic mechanical principle, i.e conservation of linear and angular momentum, rather different structures of the peridynamic kernel are employed, so that it appears surely intriguing to analyze both main shared features and differences between such different models.
\section{Features of kernel structure}
\label{sec: features}

\subsection{Nonlocality}

\label{subsec:non-loc}
As stated earlier, the key feature of peridynamics is nonlocality. In the previous section, the role of nonlocality has been highlighted by the fact that internal forces acting at a point $ \x $ at every time $ t $ are determined by a process of integration over a \textit{finite} region; but differences between the classical (local) theory of continua and peridynamics may be emphasized within the framework of a more fundamental subject: kinematics.
\begin{figure}
\centering
\includegraphics[scale=0.4]{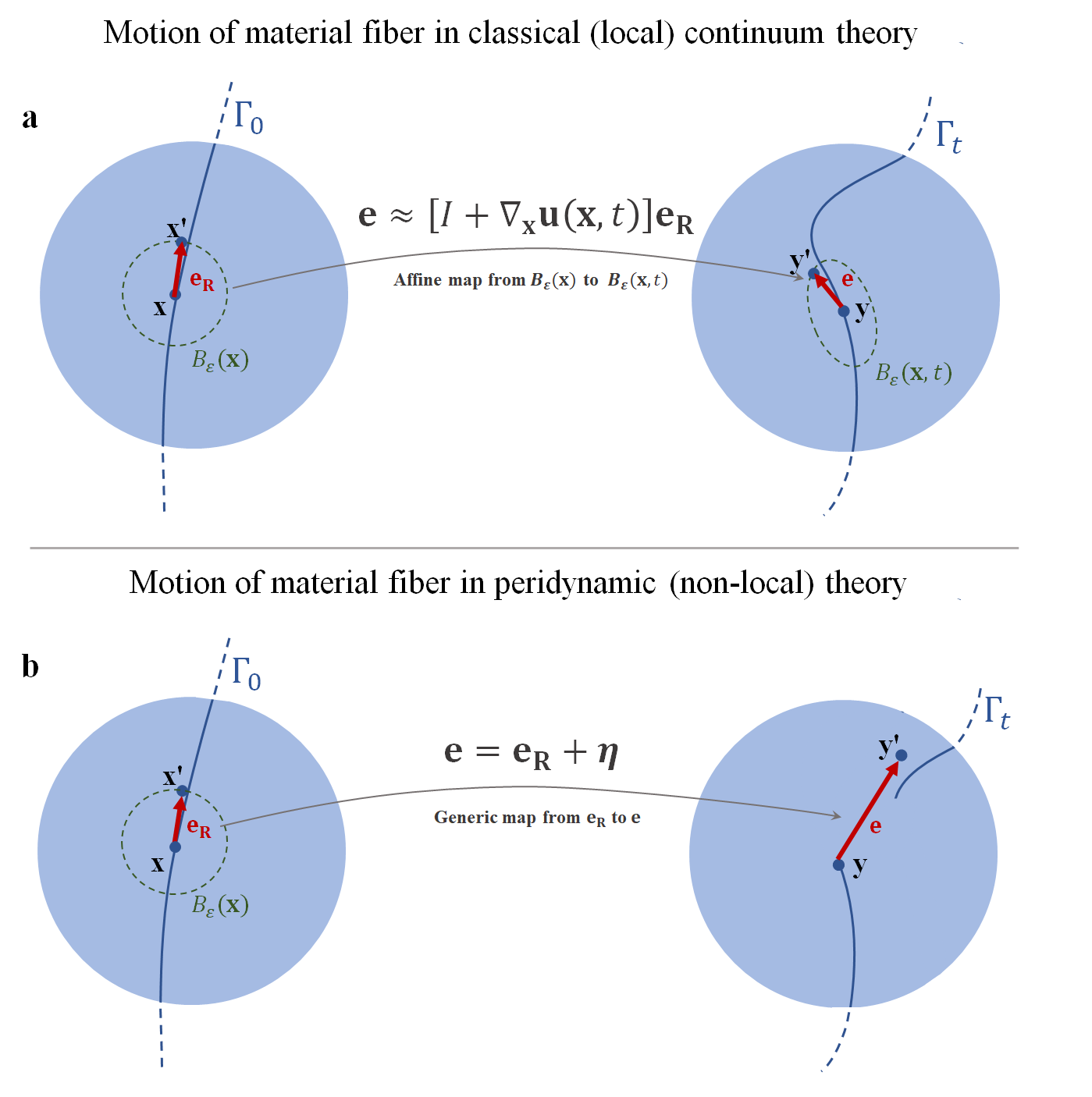}
\caption{{\bf Deformation of an infinitesimal fiber $\bm{e_R}$ connecting two points in local and nonlocal kinematics}. A material curve passing trough $ \x $, $ \Gamma_0 $, and its time evolution, $ \Gamma_t $, are represented to highlight the differences between local and nonlocal kinematics. \textbf{a)  } in classical continuum mechanics the motion of fiber is described by an affine mapping $ \bm{e}\approx [(\bm{I}+\nabla_\x\u(\x)]\bm{e_R} $ from material to geometric space; the continuous curve $ \Gamma_0 $ can never break-up. \textbf{b)} In peridynamics (nonlocal theory ) the fiber displacement is finite and determined by the generic mapping $ \bm{e}=\bm{e_R}+\eeta  $. No smoothness is required, then, a continuous material curve $ \Gamma_0 $ can break up.}
\label{fibra}
\end{figure}
The crucial difference between classical mechanics and peridynamics relies on the role of $ \nabla_\x \u $ on points motion. Given the displacement function $ \u(\x,t) $, we get that the deformed position of a point in $ \Omega_0 $ is described by the motion mapping $ \left(\x+\u(\cdot,t)\right) $, that maps points in reference space to points in geometrical space; the schism between classical local theory and peridynamics involves the Taylor expansion of the motion of a so-called \textit{undeformed fiber},  $ \x+\ell\bm{e_R} $ \cite{gurtin2010mechanics}, where $ \ell $ is the length of the fiber and $ \bm{e_R} $ a unitary vector in $\Omega_0$ representing fiber direction. Let $ \x' $ be a point in $\Omega_0$  individuated by a undeformed fiber, that is $ \x':= \x+\ell(\bm{e_R}) $.  After performing a Taylor expansion of the motion of the fiber with respect to the base point $ \x $, the following expression is obtained   
\begin{equation}
    \y'-\y=[\bm{I}+\nabla_\x\u(\x)](\x'-\x)+\text{o}(|\x'-\x|),
\end{equation}
where 
\begin{equation}
    [\nabla_\x\u(\x)]_{i,j}:=\frac{\p u(\x)_i}{\p x_j},
\end{equation}
and where $\y':=\y +\bm{e} $ (with $\bm{e}$ being, in general, non unitary) represents the corresponding \textit{deformed fiber} in $\Omega_t$ and $\bm{I}$ is the identity matrix.
Considering the limit case $ \ell \rightarrow 0 $, the  transformation law from undeformed fiber into deformed fiber results as \cite{gurtin2010mechanics} (see Fig. \ref{fibra}a)
\begin{equation}
   \bm{e}=\y'-\y=\ell([\bm{I}+\nabla_\x\u(\x)]\bm{e_R})+\text{o}(\ell)\approx \ell([(\bm{I}+\nabla_\x\u)\x]\bm{e_R}),
\end{equation}
where time dependence of $ \u $ has been omitted; transformation laws for surfaces and volumes in $\Omega_0$ are obtained similarly. Given a unitary length of undeformed fiber ($\ell=1$), so that
\begin{equation}
    \bm{e}\approx [(\bm{I}+\nabla_\x\u(\x)]\bm{e_R}
\end{equation}
holds, its stretching can be defined, as usually, by $ s:= |\bm{e}|/|\bm{e_R}|$, therefore
\begin{equation}
\begin{split}
    s:=\frac{|\bm{e}|}{\left|\bm{e_{R}}\right|}=\frac{\bm{e} \cdot \bm{e}}{1} & =(\bm{e_R}+\nabla_\x\u \bm{e_R}) \cdot(\bm{e_R}+\nabla_\x \u \bm{e_R}) \\ & =\underbrace{1+\nabla_\x\u \bm{e_R} \cdot \bm{e_R}+\nabla_\x\u^{\top} \bm{e_R} \cdot \bm{e_R}}_{\text {linear stretch}}+\underbrace{\nabla_\x \u^{\top} \nabla_\x\u \bm{e_R} \cdot \bm{e_R}}_{\text {nonlinear stretch}}.
\end{split}
\label{stretch}
\end{equation}
\textcolor{blue}{It could be noticed that the stretch of the fiber is clearly governed only by local measures. Since elastic forces are in direct relationship with the stretch, it follows that the dynamical evolutionary problem is governed by local concepts.} Therefore, in classical continuum mechanics, material bodies kinematics relies on assuming the existence of $ \nabla_\x\u $ on $ \Omega_t $ for every time $ t \in [0,+\infty) $, so that a Taylor expansion for the motion at every point in the body can be performed. In such a way, kinematics of $ \x $ and its neighborhood are defined by an affine mapping between reference and spatial points, defined by $ \nabla_\x\u $. In short, the limits of classical mechanics stem directly from the requirement that the motion of a material body must be described by a continuous time sequence of differentiable manifolds. 
Looking at peridynamics elasticity theory, it is prominent that $ \nabla_\x \u $ does not play any role in the description of points motion, so that the stretch of a fiber is purely described by finite (nonlocal) and generally nonlinear quantity, i.e. the bond stretch $ s= |\xxi+\eeta|/|\xxi| $. Looking at the peridynamics models listed above it can be seen how all pairwise forces are functions of bonds stretch, i.e. can be expressed as follow
\begin{equation}
     \f(\xxi,\eeta)=h(| \xxi+\eeta|,|\xxi|)(\xxi+\eeta)=\Tilde{h}(s)(\xxi+\eeta) ,
\end{equation}
 so that, from a finite (nonlocal) definition of fiber stretching (see Fig.\ref{fibra}b) it follows a nonlocal dynamics for every reference point in $\Omega_0$. No prescriptions are imposed on the displacement function $ \u(\x,t) $, consequently, the motion of the material body can be described by nondifferentiable, or even discontinuous mapping (with respect to the space variable); obviously, such mathematical freedom is reflected by the variety of peridynamic kernels adopted in the literature. 
\subsection{Impenetrability of matter principle}
\label{subsec:non-pen}
From the mathematical freedom allowed in peridynamics stem integral kernels of various types, sharing, however, a common structure: each model falls into the class of rational functions of variables $ \xxi+\eeta $ and $ \xxi $ (micro-modulus function structure $ c(\xxi,\delta) $, which can be quite general, can be neglected in this classification because irrelevant, due to its structure, in the following discussion on singularities and impenetrability of matter). In the following, certain features deriving from the presence of the denominator will be discussed briefly, with respect to a fundamental physical principle: the impenetrability of matter. 
To be physically consistent, the defined peridynamic kernel $ \f(\xxi,\eeta) $ must necessarily respond to the basic principle of the impenetrability of matter, that is: two material bodies cannot occupy the same position at the same time. In classical continuum theory, the condition of impenetrability is imposed via kinematics constrain on $ J:=\text{det}(\nabla_\x \chi_{t})$, where $  \chi_t$ is the deformation function ($ \chi_t:=\x+\u(\x,t) $), such that $ J>0 $. Noting that in classical continuum mechanics $\u(\x,t)$ is assumed to be a smooth function (so that $\nabla_\x \chi$ is continuous), and considering that  $ J=1$ in the reference configuration, the kinematic constrain $ J>0 $ implies the impossibility for every volume element $ \Delta v $ belonging to $ \Omega_0 $ (which is individuated by a vector basis) to change its orientation with time. Indeed, by continuity, it would imply that at a certain $ t=\tau $, $J=0$, i.e. $ \Delta v $ collapses into a plane (or a line), so that, points in the volume have been compenetrated. The same reasoning is analogous for a surface element $ \Delta s$. Obviously, for peridynamics theory, this type of constraint has no meaning. \textcolor{blue}{From a different perspective, impenetrability condition is granted by the functional form of $\Phi(\x'-\x, \mathbf{u}(\mathbf{x'}, t)-\mathbf{u}(\mathbf{x}, t))$ (see Eq.\eqref{pot6}), or better, by the choice of its singular kernel. Hence, at any time $ t $, the energy $ E[\u](t) $ associated to $ \Omega_t $, given by
\begin{equation}
E[\mathbf{u}](t):=\frac{\left\|\mathbf{u}_t(\cdot, t)\right\|_{L^{2}\left(\Omega_t\right)}^{2}}{2}+\frac{1}{2} \int_{\Omega_t} \int_{B_{\delta}(\mathbf{x})} \Phi(\x'-\x, \mathbf{u}(\mathbf{x'}, t)-\mathbf{u}(\mathbf{x}, t)) d\mathbf{x'} d \mathbf{x},
\end{equation}
satisfies
\begin{equation}
\forall \; \x , \Tilde{\x} \in \Omega_0\, , \ \ \ \ |\x-\Tilde{\x}| \rightarrow 0 \ \ \Rightarrow \ \ E[\mathbf{u}](t) \rightarrow \infty\,.
    \label{imp}
\end{equation}
Therefore, a physical configuration where two points penetrate each other is energetically inaccessible.} 
 Looking at the expression for $ E[\mathbf{u}](t) $, and given the finiteness of points velocity $ \p_t \u  $, the energetic condition of impenetrability relies upon the structure of the elastic potential $ \Phi $ (which in turn determines the peridynamic kernel structure). A necessary condition for impenetrability is the possibility of the integral terms $ \int_{B_{\delta}(\mathbf{x})} \Phi(\xxi, \eeta) d \mathbf{x'} $ to be singular.
\begin{table}[ht]
\caption{List of adopted integral kernels and their corresponding elastic potential (the function $ \mu(s,t)$ has been omitted from the model corresponding to (\ref{brittle}) because negligible in the current context); models in (\ref{nanoplane}) and (\ref{nanofibre}) are not in the table because of readability. } 
\centering 
\begin{tabular}{c c c c} 
\hline\hline    
Eq. Reference & $\Phi(\xxi,\eeta)$ & Eq. Reference & $\f(\xxi,\eeta)$  \\ [0.5ex] 
\hline \\ [0.01ex] 
 & $ \frac{c}{2}(|\xxi+\eeta|-|\xxi|)^2 $ & (\ref{mod1}) & $ c(|\xxi+\eeta|-|\xxi|) \bm{n} $ \\ [1ex] 
(\ref{pot2}) & $ \alpha(|\xxi|)\left(|\xxi+\eeta|^2-|\xxi|^2\right)^2 $, & (\ref{mod2}) & $ a(|\xxi|)\left(|\xxi+\eeta|^2-|\xxi|^2\right)(\xxi+\eeta) $  \\ [1.5ex]
 & $ c(\xxi,\delta)\;\frac{(|\xxi+\eeta|-|\xxi|)^2}{2|\xxi|} $ &(\ref{brittle}) & $ c(\xxi,\delta)\;\frac{|\xxi+\eeta|-|\xxi|}{|\xxi|}\bm{n} $  \\ [1.5ex]
 & $  c(\xxi,\delta) \frac{(|\xxi+\eeta|-|\xxi|)^2}{ 2|\xxi|^2} $ & (\ref{mod4}) & $  c(\xxi,\delta) \frac{|\xxi+\eeta|-|\xxi|}{ |\xxi|^2}\bm{n} $  \\ [1.5ex]
 & $ \frac{C(\xxi)}{r+1}(\xxi+\eeta)^{r+1} $ &(\ref{mod5}) & $ C(\xxi)(\xxi+\eeta)^r $  \\[1.5ex]
(\ref{pot6}) & $ \kappa \frac{|\xxi+\eeta|^p}{|\xxi|^{N+\alpha p}}+\Psi(\xxi,\eeta) $ & (\ref{mod6}) & $ \kappa p \frac{|\boldsymbol{\eta}+\xxi|^{p-2} }{|\boldsymbol{\xi}|^{N+\alpha p}}(\boldsymbol{\eta+\xi})+\psi(\boldsymbol{\xi}, \boldsymbol{\eta}) $  \\ [1ex] 
\hline 
\end{tabular}
\label{table:list} 
\end{table}
Looking at tab.\ref{table:list} it can be noticed that the only source of singularity is the presence of denominator of type $ |\xxi|^\beta $, for some real number $ \beta >0$; in fact, $ |\xxi|^\beta \rightarrow 0 $ during the integration process over $ B_{\delta}(\x) $. Consequently, integral kernels as in (\ref{mod1}), (\ref{mod2}), (\ref{mod5}) can never reach singular value. Focusing on the integral kernel (\ref{mod6}), which is the most general model (including both linear and nonlinear material behavior) for which wellposedness of the peridynamic equation has been proven, the energetic version of the impenetrability condition seems to be suitable; in fact, solutions of this model satisfy the following energetic condition \cite{Coclite_2018}
\begin{equation}
    E[\u](t)\leq E[\u](0) \quad \text{for a.e.}\;\; t \geq 0;
\end{equation}
so that, given a configuration with finite energy, all motions that imply compenetration are naturally excluded. The possibility for the elastic potential function of reaching a singular value as a result of the integration over $ B_\delta(\x) $ is a \textit{necessary but nonsufficient} condition for impenetrability, and a closer investigation is necessary to prove that condition (\ref{imp}) holds for the considered peridynamic kernel. 
It is interesting to highlight that such a singularity condition for $ \int_{B_{\delta}(\mathbf{x})} \Phi(\xxi, \eeta) d \mathbf{x'}$ does not hold when the process of integration is carried out by numerical quadrature methods, in fact, $ |\xxi| $ is a mesh dependent quantity and it is always a finite number. The conciliation of such a discrepancy between continuum and discrete peridynamic theory, without the trivial solution of indefinitely increasing mesh density (which is obviously unfeasible), would be surely an interesting research topic.  

\subsection{Nonlinearities and fractures}
\label{subsec: nonlincrak}
What clearly distinguish up peridynamic kernels listed in Table~\ref{table:list} is linearity/nonlinearity with respect to the relative deformed position vector $ (\xxi+\eeta) $. The nonlinearity of the integral kernel is meaningful, not only for modeling the physical behavior of those materials exhibiting a nonlinear stress-strain response, or when really interesting phenomena (like solitons  \cite{silling2016solitary}, \cite{pego2019existence}) appear, but also because it may be involved in predicting the onset of fractures in a material body. Let us consider the pairwise force equation (\ref{brittle}), reported below for convenience, as representative of linear models
\begin{equation}
    \mathbf{f}(\boldsymbol{\xi}, \boldsymbol{\eta})=c(\xxi,\delta)\;\frac{|\xxi+\eeta|-|\xxi|}{|\xxi|}\mu(s,t)\bm{n} \, ,
\end{equation}
where
\begin{equation}
\mu(s,t)=\left\{\begin{array}{ll}
1 \, , & \text { if } s\left(t^{\prime}, \xxi\right)<s_{0} \, , \\
0 \, , & \text { otherwise,}
\end{array} \text { for all } 0 \leq t^{\prime} \leq t \right. ;
\end{equation}
 we will refer to it as \textit{linear kernel}; it is noteworthy that the micro-modulus function $ c(\xxi,\delta) $ is not a source of bond stretch nonlinearity because, once a bond $ (\x,\x')\equiv\bar{\xxi} $ is fixed, $ c(\bar{\xxi},\delta)$ reduces to a constant.
\begin{figure}
\centering
\includegraphics[scale=0.4]{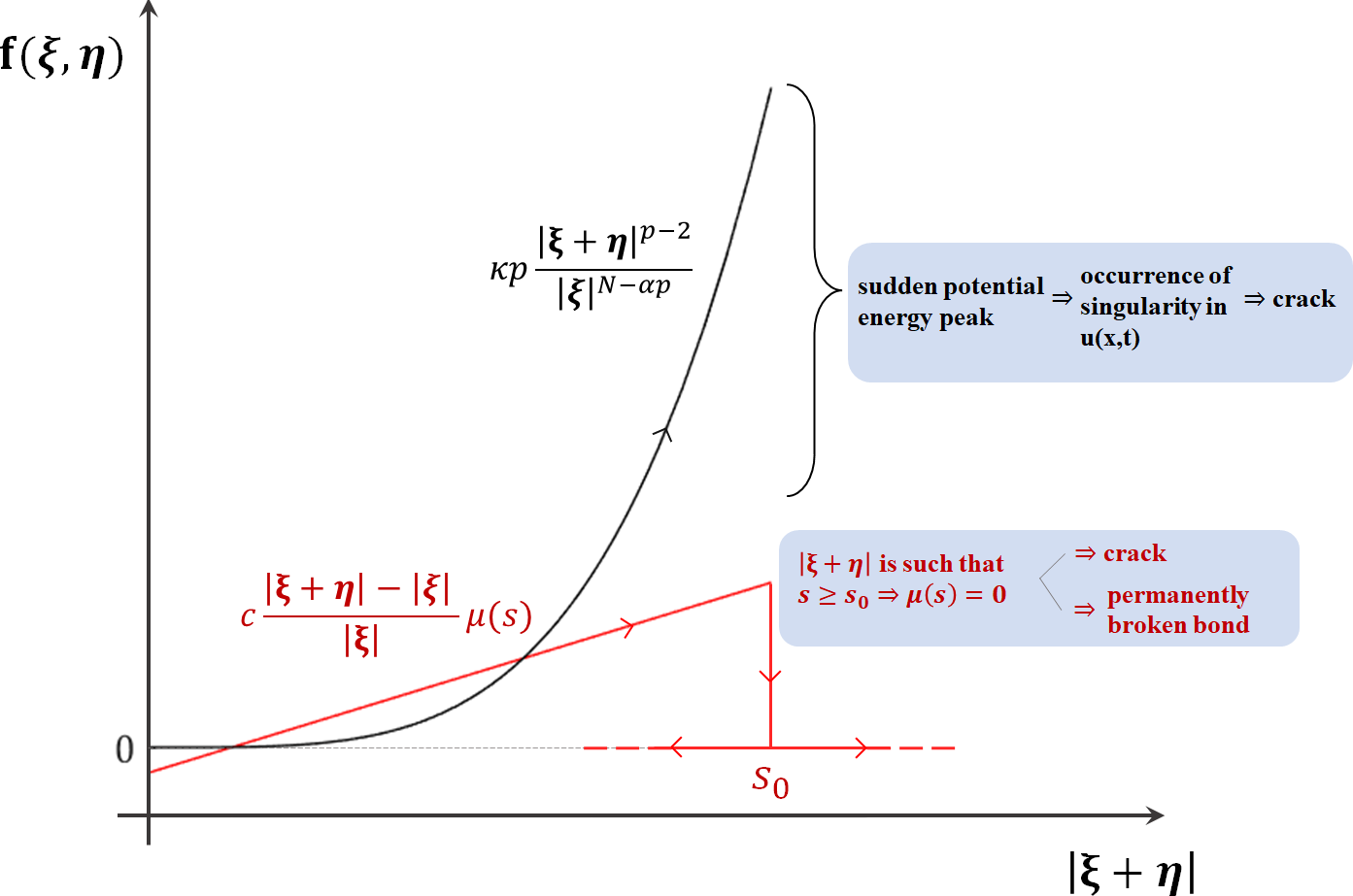}
\caption{{\bf Generation of fractures and cracks with linear and nonlinear peridynamics kernels.}  \textbf{In red} is represented a typical linear peridynamics kernel: $ c(\xxi,\delta)\;\frac{|\xxi+\eeta|-|\xxi|}{|\xxi|}\mu(s)\bm{n} $.  After a bond is broken, the pairwise force remains indefinitely equal to the zero function. \textbf{In black} is represented a nonlinear peridynamics kernel: $ \kappa p \frac{|\boldsymbol{\eta}+\xxi|^{p-2} }{|\boldsymbol{\xi}|^{N+\alpha p}}(\boldsymbol{\eta+\xi})$. Nonlinearity allows the formation of abrupt force, and so potential energy, peaks, such that, there is the occurrence of singularity in the displacement function.  }
\label{confrcrack}
\end{figure}
Let us consider, instead, as representative of nonlinear models, the pairwise force equation (\ref{mod6}), reported below for convenience 
\begin{equation}
\mathbf{f}(\boldsymbol{\xi}, \boldsymbol{\eta})=\kappa p \frac{|\boldsymbol{\eta}+\xxi|^{p-2} }{|\boldsymbol{\xi}|^{N+\alpha p}}(\boldsymbol{\eta+\xi})+\psi(\boldsymbol{\xi}, \boldsymbol{\eta});
\end{equation}
we will refer to it as \textit{nonlinear kernel}.
Regarding the linear kernel, it is notable that the only role of function $ \mu(s,t) $ is to model crack formation and development in the material, and that the threshold value $ s_0 $ is only determined by measurable quantity  \cite{silling2005meshfree}. If we consider the equivalent formulation of $ \mu(s,t) $ proposed in \cite{emmrich2016short}
\begin{equation}
\mu(s)=\Theta_{\epsilon}\left(\int_{0}^{t} \max \left\{0, s(\xxi, \tau)-s_{0}\right\} \mathrm{d} \tau\right) \quad \text{for} \; \;  \epsilon \rightarrow 0 \, ,  
\end{equation}
where
\begin{equation}
\Theta_{\varepsilon}(r)= \begin{cases}1 \, , & r \leq 0, \\ 1-\frac{r}{\varepsilon} \, , & 0<r<\varepsilon \, , \\ 0 \, , & r \geq \varepsilon \, ,\end{cases}
\end{equation}
the mathematical well-posedness of the peridynamic equation related to the linear integral kernel with bond-braking function $ \mu(s,t) $ can be proven (see \cite{emmrich2016short} for rigorous proof; for a different reformulation of the bond-braking function $ \mu(s,t) $, where no time integration is involved and for which mathematical well-posedness of the peridynamic equation is also proved, see \cite{DTT}). Notwithstanding its mathematical well-posedness, it can be noticed, looking at the linear pairwise force plot (see red line in Fig.\ref{confrcrack}), how this kind of model looks a bit artificial; in fact, a physical phenomenon described by a linear model is such that its behavior is indefinitely preserved, so that, within this setting, the onset of abrupt events (like cracks are) is unreasonable, except with a "forced from the outside" intervention.
Considering, instead, the nonlinear kernel, it can be seen how this kind of pairwise force structure allows inherently the possibility of modeling sudden force intensity (and concurrently potential energy) escalations with respect to relatively low displacements. Potential energy peaks may be naturally identified with abrupt events, like cracks, which are phenomena that require great energy localization into the bound to take place (see black line in Fig.\ref{confrcrack}). Within the context of nonlinear kernels, the aforementioned \textit{nonlinear crack} approach is capable of predicting cracks formation independently from external parameters, like $ s_0 $, and regardless of the degree of regularity of the initial condition, but, on the other hand, it is not able to propagate cracks unless a \textit{bond-breaking} function, as $ \mu(s,t) $, is not introduced.   

\section{Peridynamics in fluids: future perspectives}
\label{sec: perifluid}

Now the perspective of rigorously extending the peridynamics formalism to fluids is given. The main motivation for such an extension is to develop a model for fluids that allow to model nonlocal phenomena typically seen in complex fluid (like blood, polymeric flow, etc), and, at the same time, to take the advantage of the mathematical soundness of PD theory with respect to arising discontinuities.

As it is well-known, fluids and solids (i.e elastic materials) exhibit very different physical behavior. Let us consider a simple but clarifying example. If an elastic (solid) material is subjected to a pure shear deformation (e.g. the material is placed between two plates forced to move in opposite directions); we would expect, and it is what physically happens, that an equilibrium configuration will be reached in a certain amount of time, in the sense that the solid assume a well definite shape so as to balance plates-exerted forces. Now, if one would repeat the experiment for a fluid, as time goes on, it keeps sliding with the plates, and no equilibrium shape is reached. On the other hand, looking at the fluid velocity it can be observed that each point assumes a well-defined average velocity profile, so that an equilibrium configuration is reached in the velocity field. It can be concluded that, in solids, the "preferred" kinematic quantity is the displacement vector $\u$, while, in fluid, it is the infinitesimal variation of displacement with time, i.e the velocity $\u_t$ to assume a key role. Such a fundamental distinction reflects, in the classical continuum mechanics framework, on the structure of the governing equations for fluids and solids, respectively.

The main three ingredients in mechanical governing equations are kinematics relations, physically conserved quantities, and constitutive relations. For simplicity sake, in the following discussion, thermodynamic considerations are neglected, postponing a more comprehensive treatise for future investigations.
Kinematics of elastic materials is essentially described by the deformation gradient:
\begin{equation}
{\bf F := \nabla_{\x} \bf \chi} \, , 
\label{Kin_lag}
\end{equation}
recalling that $\chi(\x,t):= \x+\u(\x,t)$.
For the balance equations, a referential structure is adopted, so that mass, linear momentum, and angular momentum conservation read, 
\begin{equation}
\begin{gathered}
\rho=\rho_R/J, \\
\rho_{R} \ddot{\chi}=\nabla_\x \cdot {\bf T}_R+{\bf b}_{0 R}, \\
{\bf T}_R {\bf F}^{\top}={\bf F T}_{R}^{\top};
\end{gathered}
\label{BalanceLag}
\end{equation}
where $J:= \det \bf F$, and where the subscript $R$ denotes that a given quantity is considered in the reference configuration. Indeed, $\rho_R$ is the density at point $\x$ in reference configuration and $T_R$ is the first Piola stress tensor. The following constitutive relation is assumed
\begin{equation}
   { \bf T}_R=T_R(\bf F),
   \label{constit_elast}
\end{equation}
where $T_R$ is a function directly linking stress at a point of the body to its deformation gradient. Therefore, considering equations (\ref{Kin_lag}), (\ref{BalanceLag}), (\ref{constit_elast}) given below, fully describing the evolution of an elastic body, the fundamental role of $\bf F$  can be highlighted. From that, it appears clearly that $\u(x,t)$, i.e. the displacement of the material point $\x$ from its reference configuration, assumes a preferred status in elastic bodies, due to its close relation with $\bf F$, namely $\bf F := \nabla_\x \chi = \bf I + \nabla_\x \u$ .

Now let us consider kinematics, balance equations, and constitutive relation for incompressible fluids. Due to the absence of a physical meaningful reference configuration, kinematic relations are expressed only as a function of the geometrical space variables ${\bf z} \in \R^n$ (${\bf z}$ symbol instead of the classical $x$ is used to differentiate sharply the geometrical space variable from the reference space variable $\x$). Now, identifying  $\chi_t (\x,t)=\chi((\chi^{-1}({\bf z},t),t)$ with $ {\bf v }({\bf z},t)$, the motion of a given  particle of a fluid continuum is expressed as follow
\begin{equation}
    \frac{D{\bf v} ({\bf z},t)}{Dt}={\bf v}_t ({\bf z},t)+[\nabla_{{\bf z}}{\bf v} ({\bf z},t)]{\bf v} ({\bf z},t)\, ,
    \label{derivataMateriale}
   \end{equation}
where $D/Dt$ denotes the material derivative of the field ${\bf v}({\bf z},t)$ and $\nabla_{\bf z}$ is the gradient operator with respect to geometrical space variables.
The balance equations are expressed consequently \cite{gurtin2010mechanics}:
\begin{equation}
\begin{gathered}
\nabla_{\bf z}\cdot {\bf v}=0, \\
\rho {\bf v}_t=\nabla_{\bf z} \cdot ({\bf T}-\rho{\bf v}\otimes {\bf v}) + {\bf b}_0, \\
{\bf T}={\bf T}^\top;
\end{gathered} 
\label{BalanceEul}
\end{equation}
where $\bf T$ is the Cauchy stress tensor.
For fluids that are in general non-Newtonian, the following constitutive relation is assumed
\begin{equation}
    {\bf T}=T(\bf D),
    \label{cost-flu}
\end{equation}
where $T$ is a function that links the stress at a fluid point to the stretching tensor $\bf D$, defined as
\begin{equation}
    {\bf D}:= \frac{1}{2}(\nabla_{\bf z} {\bf v}+\nabla_{\bf z} {\bf v}^{\top} ),
\end{equation}
which represents the stretching and angle changing rate of material fibers. A so-called Eulerian description of fluids has been presented.
By looking at the equations that govern fluid motion, the relevance of the velocity $\bf v$ can be noticed. Being $\chi_t (\x,t)= \u_t(\x,t)$ and $\chi_t (\x,t)\equiv {\bf v}({\bf z},t)$, the key role of $\u_t$ in fluid phenomenology is recovered.
Thus, in the extension of the peridynamics formalism toward fluids, it is expected that velocity must possess a key role in the governing equation.  

A first attempt to extend PD theory to fluids could be to model the divergence of internal forces with a nonlocal term of integral type,  $\int_{B_\delta}\f dV$, similarly to what has been done for elastic materials PD theory. But, unlike what happens for the Lagrangian form of linear momentum balance equation (\ref{BalanceLag})\textsubscript{2}, where all spatial derivatives (i.e $\nabla_\x\cdot$) are replaced by the integral force term, this step would not work for the Eulerian form of linear momentum balance (\ref{BalanceEul})\textsubscript{2}, where $(\nabla_{\bf z} \cdot)$ does not appear only in the internal force term, but also is inherent in fluid kinematics. Hence, it is not surprising that proposed fluid peridynamics state-based models in \cite{silling2017modeling,behzadinasab2018peridynamics} are based on a Lagrangian description of fluid and that a truly Eulerian peridynamics formulation has not been proposed yet.

An interesting feature brings together the proposed fluid PD models: if in elastic material the peridynamics horizon $B_\delta (\x)$ is defined in the reference configuration and it is composed of the same material points during all body motion, in fluid, PD horizon maintains its spherical shape at every time $t$, so that different points at different time instant may belong to the horizon. The main explanation (see \cite{silling2017modeling}) of such a divergent property of $B_\delta(\x)$ relies on the fact that in fluids, where deformations are very huge, it would be nonphysical to let points $\x$ and $\x'$, such that $|\x-\x'|\leq \delta$ in $\Omega_0$, to keep interacting even if at a certain time, say $t=\tau$, $\y=\chi(\x)$ and $\y'=\chi(\x')$ in $\Omega_\tau$ are such that $|\y-\y'| \gg \delta$.
Such differences between fluid and solid in the definition of the PD horizon may be regarded in a different light, by the introduction of \textit{fading memory} concept. Fading memory in materials was firstly introduced by L. Boltzmann \cite{boltzmann1878theorie} and V.Volterra \cite{volterra1912equations}, and formalized by B. D. Colemann $\&$ W. Noll in the framework of continuum mechanics \cite{coleman1961foundations}; the main assumption is that the stress in a material body depends upon the history of its past configurations up to a certain remembered time. \textcolor{blue}{Specifically, the concept of fading memory is based on the following argument: the stress-free configuration for the elastic material cannot be recovered instantaneously as soon as the interaction has acted. Indeed, liquid materials instantly ``forgets'' their already acquired deformations since any state is an equilibrium state: this is why liquids are considered without memory. On the contrary, elastic materials do not forget the history of their deformations at all. So, their memory does not fade.} In this contest, the difference between elastic materials and fluids is related to their memory: an elastic material is such that it remembers its ``original'' configuration, i.e. possesses infinite memory, and thus tries to restore its initial shape at any time. On the opposite hand, a fluid is a material with zero memory, so that a reference or ``original'' configuration cannot be found \cite{astarita1975principles}.

In the following, we present a non-rigorous but interesting and stimulating perspective of the application of the concept of fading memory in the development of a peridynamics model for fluids. Consider a material body and let the parameter $s$, with the dimension of time, denotes the memory of the material, so that, the material can remember only up to time $t-s$. So, if $s\rightarrow \infty$, the material is perfectly elastic, while, if $s \rightarrow 0$, a fluid behavior is obtained. Be the displacement function extended in the time variable, so that $\Tilde{\u}(\x,t): \Omega_0 \in \R^n \times \R \rightarrow \R^n$, is such that
\begin{equation}
    \Tilde{\u}(\x,t):=\begin{cases}
    \u(\x,t) \, , &\quad \text{for} \; t\geq 0 \, , \\
    {\bf 0} \, , &\quad \text{for} \; t < 0 \, ,   \end{cases}
\end{equation}
and let us consider the consequent extension of $\chi(\x,t)$, so that $\chi(\x,t):=\x+\Tilde{\u}(\x,t)$.
 In standard PD theory, the internal force acting at point $\x$ at time $t=\tau$ is of elastic nature and directly depends on the deformation of $B_\delta(\x)$ at $t=\tau$. A straightforward extension of this peridynamics concept of internal force to materials with memory $s$, is to consider the deformation of the peridynamics horizon with respect to its last remembered configuration. i.e $B_\delta(\x,t-s)$, where 
\begin{equation}
    B_\delta(\x,t-s):=\{{\bf z} \in \R^n: \norm{\chi(\x,t-s)-{\bf z}} < \delta \}.
\end{equation}
Consequently, internal forces acting at $\x$ are related to the current deformation (at time $t$) of $B_\delta(\x,t-s)$, which is, thus, taken as a reference configuration.
It is noteworthy that, in this new framework, the horizon assumes a more geometrical structure, in contrast with the classical definition, by which, it is defined as a set of fixed material points.
At the moment, boundary are neglected, so that, the horizon is never empty, i.e $\forall \; {\bf z} \in B_\delta(\x,t) \; \exists \; \x' \in \Omega_0 \; \text{such that} \; {\bf z}=\chi(\x',t)$. 
Given the above considerations, it seems reasonable to extend the classical PD bond-based model, to a generic material with memory $s$, maintaining the same physical structure. Therefore, the pairwise force $\f$ is regarded as a function of the relative position between $\x$ and a generic point in $B_\delta(\x,t-s)$ in the reference configuration (time $t-s$), and of their current relative displacement (time $t$), with respect to the reference configuration. As a consequence, the peridynamics balance of the linear momentum equation follows:

\begin{equation}
{\begin{split}
\rho \chi_{tt}(\x, t)=\int_{\Omega_0 \cap B_{\delta}(\x,t-s)}  &\f \left( \underbrace{{\bf z}-\chi(\x,t-s)\, ,}_\text{relative position between ${\bf x}$ and ${\bf z}$ }\right. \underbrace{\chi(\chi^{-1}({\bf z},t-s),t)-\chi(\chi^{-1}({\bf z},t-s),t-s)}_\text{displacement of ${\bf z}$ from reference config. }-\\
& \left.\underbrace{ [\chi(\x,t)-\chi(\x,t-s)]}_\text{displacement of $\x$ from reference config.} \right) dV_{{\bf z}}\, .
\end{split}}
\label{memoryperi}
\end{equation}

If a solid material is considered, i.e $s\rightarrow \infty$, by assuming necessary smoothness  with respect to the time variable for all involved functions, we get 
\begin{equation}
\begin{split}
\lim_{s\rightarrow\infty}\int_{\Omega_0 \cap B_{\delta}(\x,t-s)}  &\f\left( {\bf z}-\chi(\x,t-s)\right., \chi(\chi^{-1}({\bf z},t-s),t)- \\ & \chi(\chi^{-1}({\bf z},t-s),t-s)- \left. [\chi(\x,t)-\chi(\x,t-s)] \right) dV_{{\bf z}}\\
=\int_{\Omega_0 \cap B_{\delta}(\x)}  &\f\left( \x'+\x , \;
 \chi(\x',t)-\x'- \chi(\x,t)+\x] \right) dV_{\x'}\\
 =\int_{\Omega_0 \cap B_{\delta}(\x)}  &\f\left( \x'+\x , \;
 \x'+\Tilde{\u}(\x',t)-\x'- \x -\Tilde{\u}(\x,t)+\x \right) dV_{\x'}\\
 =\int_{\Omega_0 \cap B_{\delta}(\x)}  &\f\left( \xxi , \;
 \eeta \right) dV_{\x'}\, ;
\end{split}
\label{Solidmemoryperi}
\end{equation}
remembering that, as $ (t-s) \rightarrow - \infty$, $\Tilde{\u}(\;\cdot \;, t-s) \rightarrow \bf 0 $, and, thus,  $B_\delta(\x,t-s) \ni {\bf z} \equiv \x' \in B_\delta (\x)$. The standard bond-based PD equation has been recovered.
If the limit of (\ref{memoryperi}) for $s\rightarrow0$ is taken, it is expected to obtain, almost formally, a peridynamic equation for incompressible fluids. 
Lets first consider a first-order Taylor expansion in the time variable, so that
\begin{equation}
\begin{split}
\rho\chi_{tt}(\x,t)= \int_{\Omega_0 \cap B_{\delta}(\x,t-s)} \f\left( {\bf z}-\chi(\x,t-s), \left[\chi_t(\chi^{-1}({\bf z},t-s),t)- \chi_t(\x,t) \right ]s +o(s) \right) dV_{{\bf z}}.
\end{split}
\end{equation}
Thus, under suitable regularity assumptions with respect to the time variable, letting $s\rightarrow 0$, it follows
\begin{equation}
\begin{split}
\rho\chi_{tt}(\x,t)= \int_{\Omega_0 \cap B_{\delta}(\x,t)} \f\left( {\bf z}-\chi(\x,t), \left[{ \bf v}({\bf z},t)- \chi_t(\x,t) \right ]\delta s \right) dV_{{\bf z}},
\end{split}
\label{perifluid}
\end{equation}
where $\delta s$ is an infinitesimal time interval ($\delta s \ll 1 $) we can get rid of by a slightly modification of (\ref{memoryperi}), and where the correspondence ${\bf v}({\bf z},t) \equiv  \chi_{t}(\chi^{-1}({\bf z},t))$ has been employed. 
It is interesting to note that, similarly to classical constitutive relations for fluid (\ref{cost-flu}), the actual velocity difference between $\x$ and others surrounding point in its horizon is involved in the peridynamic kernel and, thus, determine the constitutive behavior of the fluid. By employing the concept of fading memory, the key role of velocity in fluid material has been recovered, almost formally. On the other hand, a redefinition of the peridynamic horizon from a set of material points to a purely geometrical one brings extra complications. Indeed in (\ref{perifluid}) geometric and material velocity is mingled so that the model is neither Lagrangian nor Eulerian. \textcolor{blue}{It must be stressed that, the present formulation is only suitable for describing isotropic and fully-periodic flows. Future studies will include the description of the boundaries and the role of the choice of the specific strategy for enforcing no-slip conditions in the peridynamic realm.}

\section{Conclusions}

A wide-breath review of bond-based peridynamic theory has been presented. Features, applications and numerical methods involved in bond-based and state-based PD have been briefly described; then, the paper focuses on the bond-based PD formulation of solid mechanics. Basic physical principles that shape peridynamics integral kernels are described, then, a review of the most common kernels employed in the literature is displayed. By a comparison of the various kernels, it has been shown that, apart from the micro-modulus function, which can be quite general, all the models share a rational function structure. It has been shown that PD kernels derive from exact kinematics, i.e. no approximations of points displacement are assumed, at all scales. Implications of material impenetrability constrain on the kernel structures have been analyzed. A necessary, but non sufficient condition for material impenetrability has been related to the presence of a denominator of type $|\xxi|^\beta$, with $\beta$ some real number, that goes to zero in the analytical integration process. On the other hand, in a numerical integration process, this condition has no meaning (the denominator is mesh dependent and thus never reaches a zero value) so that future research in numerical peridynamics may be directed toward the resolution of such discrepancy. The potential role of employing a nonlinear PD kernel in predicting the onset of fractures has been explored. Cracks formation modeling capacity of widely employed linear kernels supplemented with a bond-breaking function is qualitatively compared with a nonlinear kernel. As a result, it appears that the presence of nonlinearity in the relative displacement can reproduce inherently the high potential energy peaks that are indicative of cracks formation process. Further numerical studies on nonlinear kernels will be performed to test their accuracy in predicting crack formation. Finally, some research perspectives on the extension of bond-based PD toward fluids have been offered. The prominent behavior difference between solids and fluids has been related to the preferred status of deformation, in solids, and deformation rate, in fluids. Then, the implications of extending the peridynamics formalism toward fluids have been analyzed via the materials memory concept.  
In this unifying framework, the main ideas of peridynamics have been extended formally toward fluids. The key role of the deformation rate has been recovered, but, on the other end, the necessary redefinition of the PD horizon from a material set to a geometric one leads to a PD governing equation in which geometric and material quantities are blended. Future research will be directed toward a more clear understanding of the implications of PD nonlocality on the relation between Lagrangian (material) and Eulerian (geometric) descriptions of materials. 

\section{Declaration statement}

\subsection*{Ethics approval and consent to participate.} The authors approve the ethics of the journal and give the consent to participate.

\subsection*{Consent for publication.} The authors consent for the publication.

\subsection*{Availability of data and material.} The datasets used and/or analyzed during the current study are available from the corresponding author on reasonable request.

\subsection*{Competing interests.} The authors declare no competing interests.

\subsection*{Funding.} This work was partially supported by the project ``Research for Innovation'' (REFIN) - POR Puglia FESR FSE 2014-2020 - Asse X - Azione 10.4 (Grant No. CUP - D94120001410008) and Fundacao para a Ciencia e a Tecnologia (FCT) under the program “Stimulus” with grant no. CEECIND/04399/2017/CP1387/CT0026 and through the research project with ref. number PTDC/FIS-AST/0054/2021.

\subsection*{Authors' contributions.} The authors declare that they gave their individual contributions in every section of the manuscript. All authors read and approved the final manuscript.

\subsection*{Acknowledgments} AC acknowledges the project ``Research for Innovation'' (REFIN) - POR Puglia FESR FSE 2014-2020 - Asse X - Azione 10.4 (Grant No. CUP - D94120001410008). AC and TP are members of Gruppo Nazionale per il Calcolo Scientifico (GNCS) of Istituto Nazionale di Alta Matematica (INdAM). GF acknowledges the support by Fundacao para a Ciencia e a Tecnologia (FCT) under the program “Stimulus” with grant no. CEECIND/04399/2017/CP1387/CT0026 and through the research project with ref. number PTDC/FIS-AST/0054/2021. GF is also a member of the Gruppo Nazionale per la Fisica Matematica (GNFM) of the Istituto Nazionale di Alta Matematica (INdAM).


\end{document}